\documentclass[]{article}

\title{Lifting a $5$-dimensional representation of $M_{11}$ to a complex unitary representation of a certain amalgam}
\author{Geoffrey R. Robinson}
\usepackage{latexsym,amssymb}
\begin{document}

\maketitle

\begin{abstract}
We lift the $5$-dimensional representation of $M_{11}$ in characteristic $3$ to a unitary complex representation of the amalgam ${\rm GL}(2,3)*_{D_{8}}S_{4}$.
\end{abstract}

\section{The representation}
It is well known that the Mathieu group $M_{11},$ the smallest sporadic simple group, has a $5$-dimensional (absolutely) irreducible representation over ${\rm GF}(3)$ (in fact, there are two mutually dual such representations). It is clear that this does not lift to a complex representation, as $M_{11}$ has no faithful complex character of degree less than $10.$ 

\medskip
However, $M_{11}$ is a homomorphic image of the amalgam $G = {\rm GL}(2,3)*_{D_{8}}S_{4},$ and it turns out that if we consider the $5$-dimension representation of $M_{11}$ as a representation of $G$, then we may lift that representation of $G$ to a complex representation. We aim to do that in such a way that the lifted representation is unitary, and we realise it over $\mathbb{Z}[\frac{1}{\sqrt{-2}}],$ so that the complex representation admits reduction (mod $p$) for each odd prime. These requirements are stringent enough to allow us explicitly exhibit representing matrices. It turns out that reduction (mod $p$) for any odd prime $p$ other than $3$ yields either a $5$-dimensional special linear group or a $5$-dimensional special unitary group, so it is only the behaviour at the prime $3$ which is exceptional.

\medskip
We are unsure at present whether the $5$-dimensional complex representation of $G$ is faithful (though it does have free kernel), so we will denote the image of $G$ in ${\rm SU}(5,\mathbb{Z}[\frac{1}{\sqrt{-2}}] )$ by $L$,
and denote the image of $L$ under reduction (mod $p$) by $L_{p}.$ 

\medskip
We recall that to construct a $5$-dimensional representation of $G$, we need to construct $5$-dimensional representations of $H = {\rm GL}(2,3)$ and $K = S_{4}$ which agree on a common dihedral subgroup of order $8$.

\medskip
We recall that $H$ has a presentation: $$\langle b,c : b^{2} = c^{3} = (bc)^{8}= [b,(bc)^{4}] = [c,(bc)^{4}] = 1 \rangle,$$ for this is a presentation of a double cover of $S_{4}$ in which the pre-image of a transposition has order $2$. It is also helpful in what follows to note that a unitary $2 \times 2$ matrix of trace $\pm \sqrt{-2}$ and determinant $-1$ has order $8$ and that a unitary $2 \times 2$ matrix of trace $-1$ and determinant $1$ has order $3$. We set $a = \left(\begin{array}{clcrc} -1& 0 & 0 & 0 & 0\\0&-1&0&0&0\\0&0&1&0&0 \\0& 0 &0 &0 & -1\\0&0&0&1& 0 \end{array} \right)$, $b = \left(\begin{array}{clcrc} -1& 0 & 0 & 0 & 0\\0&1&0&0&0\\0&0&1&0&0 \\0& 0 &0 &1 & 0\\0&0&0&0& -1 \end{array} \right)$,\\
$c = \left(\begin{array}{clcrc} \frac{1}{2}& \frac{-1}{2} & \frac{-1}{\sqrt{-2}} & 0 & 0\\ \frac{1}{2}& \frac{-1}{2} & \frac{1}{\sqrt{-2}} & 0 & 0\\ \frac{1}{\sqrt{-2}}& \frac{1}{\sqrt{-2}} & 0 & 0 & 0\\0& 0 &0 & \frac{-1-\sqrt{-2}}{2}& \frac{-1}{2}\\0&0&0&\frac{1}{2}& \frac{-1+\sqrt{-2}}{2} \end{array} \right)$,
$d = \left(\begin{array}{clcrc} 0 & 0 & 0 & 1 & 0\\0 & \frac{-1}{2} & \frac{-1-\sqrt{-2}}{2}&0&0 \\0& \frac{1-\sqrt{-2}}{2} & \frac{-1}{2}&0&0\\ 0 & 0 & 0 & 0 & 1\\1&0&0&0&0 \end{array} \right).$\\

\medskip
We note that $a$ has order $4$, that $b$ has order $2$, and that $c$ and $d$ each have order $3$. Also, $bc$ has order $8$, and $(bc)^{4}$ commutes with both $b$ and $c$. Hence $H = \langle b,c \rangle \cong {\rm GL}(2,3).$

\medskip
It is clear that $K = \langle a,b,d \rangle \cong S_{4},$ since $ad$ has order $2$. Also, $a = (c^{-1}bc^{-1})^{2}$. Hence $H \cap K \geq \langle a,b \rangle$. But $K \not \subseteq H,$ since there are $H$-invariant subspaces which are not $K$-invariant. Hence $H \cap K = \langle a,b \rangle$ is dihedral of order $8,$ so $L = \langle H,K \rangle$ is a homomorphic image of $G$ via this representation. Furthermore, the kernel of the homomorphism is free as ${\rm GL}(2,3)$ and $S_{4}$ are faithfully represented. Note that, although the generator $a$ is redundant, (as is the generator $b$), the presence of $a$ and $b$ makes it clear that $L$ is a homomorphic image of the amalgam $G$. 

\section{Reductions (mod $p$)}

\medskip
We now discuss the groups $L_{p}$, where $p$ is an odd prime. More precisely, we reduce the given representation (mod $\pi$), where $\pi$ is a prime ideal of $\mathbb{Z}[\sqrt{-2}]$ containing the odd rational prime $p.$
It is clear that $L_{3}$ is a subgroup of ${\rm SL}(5,3)$ (and choosing different prime ideals containing $3$ leads to representations dual to each other). Computer calculations with GAP confirm that $L_{3} \cong M_{11}.$
(I am indebted to M. Geck for assistance with this computation). Suppose from now on that $p >3.$ If $p \equiv 1$ or $3$, (mod $8$), then $-2$ is a square in ${\rm GF}(p).$ If $p \equiv 5$ or $7$, (mod $8$), then $-2$ is a non-square in ${\rm GF}(p).$ Hence $L_{p}$ is a subgroup of ${\rm SL}(5,p)$ when $p \equiv 1$ or $3$ (mod $8$) and $L_{p}$ is a subgroup of ${\rm SU}(5,p)$ when $p \equiv 5$ or $7$ (mod $8$). We will prove:

\newpage
\noindent {\bf Theorem 1}

\medskip
\noindent \emph{i) $L_{3} \cong M_{11}$}\\
\noindent \emph{ii)$L_{p} \cong {\rm SL}(5,p)$  when $p > 3$ and $p \equiv 1$ or $3$ (mod $8$).}\\
 \noindent \emph{iii) $L_{p} \cong {\rm SU}(5,p)$ when $p \equiv 5$ or $7$ (mod $8$).}

\medskip
\noindent {\bf Remarks:} We note, in particular, that the Theorem implies that $L$ is infinite, although we need to establish this fact during the proof in any case. We also note that $G$ is not isomorphic to ${\rm SU}(5, \mathbb{Z}[\frac{1}{\sqrt{-2}}]),$ since $G$ contains no elementary Abelian subgroup of order $8$ (since it is an amalgam of finite groups, neither of which contains such a subgroup), but ${\rm SU}(5,\mathbb{Z}[\frac{1}{\sqrt{-2}}])$ contains elementary Abelian subgroups of order $16$. In fact, the theorem also implies that $L$ is not isomorphic to ${\rm SU}(5,\mathbb{Z}[\frac{1}{\sqrt{-2}}]),$ since all elementary Abelian $2$-subgroups of $L$ map isomorphically into $L_{3}$, and $L_{3}$ contains no elementary Abelian subgroup of order $8$. We recall, however, that, as noted in [5],  J-P. Serre has proved that $G$ is isomorphic to ${\rm SU}(3, \mathbb{Z}[\frac{1}{\sqrt{-2}}]).$ 

\medskip
We note also that $G$ has the property that all of its proper normal subgroups are free. Otherwise, there is such a normal subgroup $N$ that contains an element of order $2$ or an element of order $3$. All involutions in $G$ are conjugate, because $G$ has a semi-dihedral Sylow $2$-subgroup with maximal fusion system. Both $S_{4}$ and ${\rm GL}(2,3)$ are generated by involutions so if $N$ contains an involution, we obtain $N = G.$ Now $G$ has two conjugacy classes of subgroups of order $3,$ so if $N$ contains an element of order $3,$ then $N$ contains a subgroup isomorphic to $A_{4}$ or to ${\rm SL}(2,3)$, so contains an involution, and $N = G$ in that case too.

\medskip
Now we proceed to prove that $L$ is infinite. It is clear that $L$ is irreducible, and primitive, as a linear group. We will prove more generally that no finite homomorphic image of $G$ has a faithful complex irreducible representation of degree $5$. If $M$ were such a homomorphic image then we would have $M = [M,M]$ and $M$ is primitive as a linear group (otherwise $M$ would have a homomorphic image isomorphic to a transitive subgroup of $S_{5}$, which must be isomorphic to $A_{5,}$ as $M$ is perfect). But $M \cong G/N$ for some free normal subgroup $N$ of $G,$ so that $M$ has subgroups isomorphic to $S_{4}$ and ${\rm GL}(2,3),$ a contradiction.

\medskip
Now R. Brauer (in [2]), has classified the finite primitive subgroups of ${\rm GL}(5,\mathbb{C}),$ so we make use of his results. If $O_{5}(M) \not \subseteq Z(M)$, then $M/O_{5}(M)$, being perfect, must be isomorphic to ${\rm SL}(2,5)$, since $O_{5}(M)$ is irreducible, and has a critical subgroup of class $2$ and exponent $5$ on which elements of $M$ of order prime to $5$ act non-trivially. But $M/O_{5}(M)$ contains  an isomorphic copy of ${\rm GL}(2,3),$ a contradiction, as ${\rm SL}(2,5)$ has no element of order $8$.

\medskip
Hence $M$ must be isomorphic to one of $A_{6},{\rm PSU}(4,2)$ or ${\rm PSL}(2,11)$. We have made use of the fact that the $5$-dimensional irreducible representation of $A_{5}$ is imprimitive. We also use transfer to conclude that $Z(M)$ is trivial. Since $M = [M,M],$ we see that the given representation is unimodular, so $Z(M)$ has order dividing $5$. But since $M/Z(M)$ has a Sylow $5$-subgroup of order $5,$ when $S$ is a Sylow $5$-subgroup of $G,$ we have $Z(M) \cap S = M^{\prime} \cap Z(M) \cap S \leq S^{\prime} = 1,$ as $S$ is Abelian. Now none of 
$A_{6},{\rm PSU}(4,2)$ or ${\rm PSL}(2,11)$ contain an element of order $8$, whereas $M$ contains a subgroup isomorphic to ${\rm GL}(2,3),$ and does contain an element of order $8$. Hence $M$ must be infinite, as claimed
(we note that Brauer's list contains $O_{5}(3)^{\prime},$ but this is isomorphic to ${\rm PSU}(4,2),$ which we have dealt with, and the realization as ${\rm PSU}(4,2)$ makes it clear that it can contain no element of order $8$).

\medskip
Now we proceed to prove that $L_{p}$ is as claimed for primes $p > 3.$  We note that $L_{p}$ has order divisible by $p$ since otherwise $L_{p}$ is isomorphic to a finite subgroup of ${\rm GL}(5,\mathbb{C}),$ which we have excluded above, as $L_{p}$ is a homomorphic image of $G.$ Now $L_{p}$ is clearly absolutely irreducible as a linear group in characteristic $p$, and $L_{p}$ is also primitive as a linear group, since we have already noted that no homomorphic image of $G$ is isomorphic to a transitive subgroup of $S_{5}.$ Let $F_{p}$ denote the Fitting subgroup of $L_{p}$. If $F_{p}$ is not central in $L_{p}$, then $F_{p}$ must be a non-Abelian $5$-group, and we see that $L_{p}/F_{p}$ is isomorphic to ${\rm SL}(2,5),$ a contradiction, as before. Thus $L_{p}$ has a component $E_{p} = E,$ which still acts absolutely irreducibly by Clifford's Theorem. Hence the component $E$ is unique. Since $L_{p}$ is perfect, and $L_{p}/E$ is solvable (using the Schreier conjecture), we see that $E = L_{p},$ and that $L_{p}$ is quasi-simple. It is clear that $L_{p}$ is a subgroup of ${\rm SL}(5,p)$ if $p \equiv 1,3$ (mod $8$), and a subgroup of ${\rm SU}(5,p)$ if $p \equiv 5,7$ (mod $8$).

\medskip
By a slight abuse, we still let $a,b,c,d$ denote their images in $E,$ for ease of notation. We note that $X  = C_{E}(a^{2})$ is still completely reducible, since it acts irreducibly on each eigenspace of $a^{2}.$ Hence $O_{p}(X) = 1.$ Suppose that $X$ contains an element $y$ of order $p.$ Then since $p \geq 5,$ $y$ must centralize $F(X)$ by the Hall-Higman Theorem. Since $O_{p}(X) = 1$ , $X$ must have a component, $T,$ say.
If $T$ has a unique involution, say $t$, then $t$ acts trivially on the $1$-eigenspace of $a^{2}$ by unimodularity, so $t$ must act as multiplication by $-1$ on the $-1$ eigenspace of $a^{2},$ and in fact $t = a^{2}$. Furthermore, $T$ must act faithfully on the $-1$-eigenspace of $a^{2}$, so that $T \cong {\rm SL}(2,p)$
in that case.

\medskip
Suppose that $L_{p}$ contains no elementary Abelian subgroup of order $8$. Then results of Alperin, Brauer and  Gorenstein ([1]) show that $L_{p}$ is isomorphic to an odd central extension of  $M_{11}, {\rm PSU}(3,q)$, or ${\rm PSL}(3,q)$ for some odd $q$. We have excluded groups with a Sylow $2$-subgroup isomorphic to a Sylow $2$-subgroup of ${\rm PSU}(3,4)$ since $L_{p}$ contains elements of order $8.$ Also, we know that $L_{p}$ contains a semi-dihedral subgroup of order $16$, so $L_{p}$ does not have a dihedral Sylow $2$-subgroup. Note also that $L_{p}$ has centre of order dividing $5$ by unimodularity. We note that since $L_{p}$ contains elements of order $p,$ we can only have $L_{p} \cong M_{11}$ if $p =5$ or $11$ (and in that case, $L_{p}$ has trivial centre by a transfer argument). In fact,  using [3], for example, $M_{11}$ has no faithful $5$-dimensional representation in any characteristic other than $3$, so we can exclude that possibility. Likewise, we do not need to concern ourselves with ${\rm PSL}(3,3)$ or ${\rm PSU}(3,3)$, using the Modular Atlas ([3]). In the other cases, every involution of $\tilde{ L}_{p} = L_{p}/Z(L_{p})$ has a component ${\rm SL}(2,q)$ (note that $\tilde{L}_{p}$ has a single conjugacy class of involutions). In fact, it follows from inspection of the given representation that every involution of $L_{p}$ has a component isomorphic to ${\rm SL}(2,q)$, since a central element of order $5$ does not have unimodular action on any eigenspace of an involution. Now let $q = r^{m}$ for some odd prime $r$. If $r \neq p,$  then ${\rm SL}(2,r)$ has a $2$-dimensional complex representation so $r \leq 5.$ However, we can exclude $r \leq 5 $ using [3]. This leaves $r =p,$ and $\tilde{L}_{p} \cong {\rm PSL}(3,p)$ or ${\rm PSU}(3,p).$ However, for $p >5,$ as noted by R. Steinberg, the Schur multiplier of ${\rm PSL}(3,p)$ or ${\rm PSU}(3,p)$ has order dividing $3$, and (using [4], for example), the only non-trivial irreducible modules of dimension less than $6$ for either of these groups are the natural module and its dual (note that the dual is also the Frobenius twist in the unitary case).

\medskip
Suppose then that $L_{p}$ contains an elementary Abelian subgroup of order $8.$ Then $L_{p}$ contains an involution $t$ which has the eigenvalue $-1$ with multiplicity $4$ and the eigenvalue $1$ with multiplicity $1$
(the Brauer character can't take the value $1$ on every non-identity element of an elementary Abelian subgroup of order $8$). Then $L_{p} \times \langle -I \rangle $ is generated by its reflections.

\medskip
By the results of Zalesskii and Serezhkin [6], we may conclude that $L_{p} \cong {\rm SL}(5,p)$ or ${\rm SU}(5,p).$ Several of the options from [6] are eliminated in our situation. For example, we have already that $L_{p}$ is not liftable to a finite complex linear group, and it is clear that $L_{p}$ is not a covering group of an alternating group (for such an alternating group would have to be of degree at most $7$ and contains no element of order $8$). We also note that $L_{p}$ is not conjugate to an orthogonal group in odd characteristic, because $bc$ is an element of order $8$ whose eigenvalues other than $-1$ do not occur in mutually inverse pairs. Its eigenvalues are $-1,\alpha^{2},\alpha^{-2}, \alpha,\alpha^{3}$ for some primitive $8$-th root of unity $\alpha.$

\section{Concluding remarks}
One way to see that $L_{3}$ is isomorphic to $M_{11}$ is to reduce the representation modulo the ideal $( 1 + \sqrt{-2})$, which clearly realizes $L_{3}$ as a subgroup of ${\rm SL}(5,3).$ It turns out that $L_{3}$ has one orbit of length $11$ on the $1$-dimensional subspaces of the space acted upon (the other orbit being of length $110$), and the resulting permutation group on the $11$ subspaces of that orbit is $M_{11}$. In reality, it is knowledge of this representation which led to the attempt to lift it to a complex representation of the amalgam.

\medskip
As we remarked earlier, we are unsure at present whether the representation of $G$ afforded by $L$ is a faithful one. Consequently, while we know that all proper normal subgroups of $G$ are free, we have not proved that this is the case for $L$. We therefore feel it is worth noting:

\medskip
\noindent {\bf Theorem 2:} \emph{Neither $G$ nor $L$ has any non-identity solvable normal subgroup.}

\medskip
\noindent {\bf Proof:} This is clear for $G,$ but for completeness we indicate a proof. Every proper normal subgroup of $G$ is free. Hence if $1 \neq S \lhd G,$ is solvable, then $S$ is free of rank one. But $G = [G,G],$ so that $S \leq Z(G)$. Now suppose that there is a non-identity element $s \in S$, and recall that $G$ has the form $H*_{D}K,$ where $H \cong {\rm GL}(2,3)$, $K \cong S_{4}$ and $D = H \cap K$ is dihedral with $8$ elements.
Now since $s$ has infinite order, $s$ may be expressed in the form $s = d x_{1}x_{2} \ldots x_{m} x_{m+1}$, where $d \in D, m \geq 1$ and each $x_{i}\in (H \cup K) \backslash D$ but there is no value of $i$ for which both $x_{i}$ and $x_{i+1}$ both lie in $H,$ and there is no value of $i$ for which $x_{i}$ and $x_{i+1}$ both lie in $K$. The expression is not unique, but for each $i,$ the right coset of $D$ containing $x_{i}$ (in whichever of $H$ or $K$ contains $x_{i}$) is uniquely determined.

\medskip
But for any $c \in D$, we have $s = s^{c} = d^{c}x_{1}^{c}x_{2}^{c} \ldots x_{m+1}^{c}$. 
It follows that $x_{i}^{c}x_{i}^{-1} \in D$ for each $i$ and each $c \in D.$ Hence each $x_{i}$ normalizes $D$.
But $D$ is self-normalizing in $K$ and $N_{H}(D)$ is semi-dihedral of order $16,$ so that $s \in N_{H}(D),$ 
a contradiction, as $s$ has infinite order.

\medskip
As for $L,$ note that if $S \lhd L$ is solvable, then $[L,S]$ is in the kernel of each reduction (mod $p$), as $L_{p}$ is always quasi-simple. However, given a matrix $x \in L$, there is a minimal  non-negative integer $s$ such that $2^{s} x$ has all its entries in $\mathbb{Z}[\sqrt{-2}]$. Now if $x \neq I,$ then there are only finitely many prime ideals of $\mathbb{Z}[\sqrt{-2}]$ which contain all entries of $2^{s}x - 2^{s}I.$ Hence $[L,S] = I.$ But, as $L$ is an irreducible linear group, $Z(L)$ consists of scalar unitary matrices of determinant $1$ with entries in  $\mathbb{Q}[\sqrt{-2}],$ so $Z(L) = 1.$

\medskip
\noindent {\bf Remark:} It might also be worth noting that Theorem 1 implies that the only torsion that $L$ can have is $2$-torsion, $3$-torsion, or $5$-torsion. Only elements of $3$-power order can be in the kernel of reduction (mod $3$), so the only possibilities for prime orders of elements of $L$ are $2,3,5$ or $11$. But any element of order $11$ in $L$ would have trace an irrational element of $\mathbb{Q}[\sqrt{-11}]$, while its trace must be in $\mathbb{Q}[\sqrt{-2}]$. At present, we see no obvious way to prove that $L$ has no $5$-torsion, since $L_{p}$ always contains elements of order $5$. We do note that $L$ does not contain the obvious permutation matrix $f = \left(\begin{array}{clcrc} 0 & 1 & 0 & 0 & 0\\ 0&0& 1& 0 & 0 \\0&0&0&1&0\\0&0&0&0&1\\1&0&0 &0&0\end{array} \right),$ since $\langle b,f \rangle$ contains an elementary Abelian subgroup of order $16$ and $L$ does not.

\medskip
\noindent {\bf Acknowledgement:} I am indebted to J.E. Humphreys for pointing out the reference [4], as well as the existence of some related theory due to A. Premet, to me.

\medskip
\begin{center}
{\bf Bibliography}
\end{center}

\medskip
\noindent  [1] Alperin, J. L.; Brauer, Richard; Gorenstein, Daniel,\emph{ Finite simple groups of 2-rank two}, Collection of articles dedicated to the memory of Abraham Adrian Albert, Scripta Math. 29, no. 3-4, (1973),191-214.  

\medskip
\noindent [2] Brauer, Richard ,\emph{\"Uber endliche lineare Gruppen von Primzahlgrad}, 
Math. Ann. {\bf 169}, (1967), 73-96. 

\medskip
\noindent [3] Jansen, Christoph; Lux, Klaus; Parker, Richard; Wilson, Robert, \emph{An atlas of Brauer characters},  (Appendix 2 by T. Breuer and S. Norton), London Mathematical Society Monographs. New Series, {\bf 11}, Oxford Science Publications, The Clarendon Press, Oxford University Press, New York, 1995.

\medskip
\noindent [4] 
L\"ubeck, Frank,\emph{ Small degree representations of finite Chevalley groups in defining characteristic}, LMS J. Comput. Math. 4 (2001), 135-169 (electronic). 
 
\medskip
\noindent [5] Robinson, Geoffrey R., \emph{Reduction mod q of fusion system amalgams}, Trans. Amer. Math. Soc. 363, {\bf 2}, (2011), 1023-1040.

\medskip
\noindent [6] Zalesskii, A. E.; Serezhkin, V. N.,\emph{Finite linear groups generated by reflections}, Izv. Akad. Nauk SSSR Ser. Mat. 44, {\bf 6},38, (1980), 1279-1307.

\end{document}